\documentclass[sn-mathphys-num]{sn-jnl}

\usepackage{graphicx}%
\usepackage{multirow}%
\usepackage{amsmath,amssymb,amsfonts}%
\usepackage{amsthm}%
\usepackage{mathrsfs}%
\usepackage[title]{appendix}%
\usepackage{xcolor}%
\usepackage{textcomp}%
\usepackage{manyfoot}%
\usepackage{booktabs}%

\usepackage{listings}%

\usepackage{amsmath,amsfonts,amsbsy,amssymb,amsthm,mathtools,tabulary,graphicx,caption,fancyhdr, mathrsfs, pgf, pgfplots, epsfig, epstopdf, lscape,booktabs,enumitem,lipsum,tkz-graph,scalerel,stackengine,setspace,adjustbox,float,amsthm}
\usepackage[labelformat=simple]{subfig}

\theoremstyle{thmstyleone}%
\theoremstyle{thmstyletwo}%

\theoremstyle{thmstylethree}%

\begin{document}

\title{Statistical Geometry and Information Dynamics on Hyperspherical Surfaces}


\author*[1]{\fnm{Masoud} \sur{Ataei}}\email{masoud.ataei@utoronto.ca}


\affil[1]{\normalsize\orgdiv{Department of Mathematical and Computational Sciences},\\ \orgname{University of Toronto}, \state{Ontario}, \country{Canada}}




\abstract{	
We study the statistical geometry of random chords on n-dimensional spheres by deriving explicit analytical expressions for the chord length distribution and its associated structural properties. A critical threshold emerges at dimension 19, marking the transition from curvature-dominated variability to high-dimensional concentration, where interpoint distances become nearly deterministic and probabilistic diversity collapses into geometric uniformity.
We further derive a closed-form expression for the Fisher information, showing that it is inversely proportional to the square of the radius and varies non-monotonically with dimension. Notably, it attains a minimum at dimension 7, coinciding with the dimension at which the volume of the unit sphere is maximized. This reflects a unique regime of maximal spatial diffuseness and minimal inferential sensitivity, arising from the interplay between weakening curvature and still-latent concentration. The alignment between volumetric and statistical extrema reveals a deeper duality between geometry and information.
We also analyze the characteristic function, which exhibits a dichotomy: in even dimensions, it takes rational-exponential form, while in odd dimensions it involves Bessel and Struve functions. This distinction reflects differences in harmonic structure and boundary regularity across dimensions. Together, these findings show how curvature and dimension jointly regulate statistical efficiency on hyperspherical domains, with implications for geometric inference and high-dimensional learning.
}




\maketitle

\newpage

\section{Introduction}

High-dimensional geometry plays a foundational role in modern statistics, information theory, and learning theory. Among canonical geometric structures, the \( n \)-dimensional sphere occupies a central position due to its symmetry, compactness, and ubiquity in applications ranging from directional data analysis to representation learning on manifolds.

A fundamental geometric question concerns the statistical behavior of interpoint distances on the hypersphere. Specifically, if two points are sampled independently and uniformly from the surface of a sphere, what is the distribution of their Euclidean separation? This chord length distribution encapsulates both the intrinsic curvature of the space and the emergent regularity of high-dimensional geometry. As the dimension increases, these distributions exhibit a transition, shifting from curvature-induced dispersion to high-dimensional concentration around a deterministic value. Understanding this phenomenon is critical not only for geometric inference but also for the design of statistical and algorithmic frameworks that operate over hyperspherical domains.

To investigate these phenomena, we derive explicit characterizations of the statistical geometry of random chords on \( n \)-dimensional spheres. We obtain closed-form expressions for the chord length distribution and associated statistical quantities including the cumulative distribution function, raw moments, Fisher information, and characteristic function. We identify dimension \( n = 19 \) as a critical threshold marking the onset of high-dimensional concentration, and analyze the implications of this transition for probabilistic modeling. Additionally, we show that the Fisher information varies non-monotonically with dimension, attaining a minimum near \( n = 7 \), where the volume of the unit sphere is also maximized. This coincidence highlights a regime of maximal spatial diffuseness and minimal inferential sensitivity, suggesting a geometric-inferential duality.

We also reveal a structural dichotomy in the analytic form of the characteristic function: in even dimensions, it assumes a rational-exponential form, whereas in odd dimensions it involves Bessel and Struve functions. This asymmetry reflects underlying differences in harmonic decomposition and analytic regularity across dimensions.

The paper is organized as follows. In Section~\ref{Sec:CLD}, we derive key aspects of the chord length distribution, including its probability density function, cumulative distribution function, raw moments, critical dimension, Fisher information, and characteristic function. In Section~\ref{Sec:Results}, we provide formal analytic derivations of these results. Finally, Section~\ref{Sec:Conclusion} offers concluding remarks and outlines directions for future work.

\section{Chord Length Distribution}
\label{Sec:CLD}
\subsection{Definition and Properties}

Let \( (\Omega, \mathcal{F}, \mathbb{P}) \) be a probability space, where \( \Omega \) is the sample space, \( \mathcal{F} \) a \( \sigma \)-algebra of measurable subsets, and \( \mathbb{P} \) a probability measure on \( \mathcal{F} \). Consider a real-valued random variable \( X : \Omega \rightarrow \mathbb{R} \), measurable with respect to the Borel \( \sigma \)-algebra \( \mathcal{B}_{\mathbb{R}} \), with absolutely continuous distribution and associated probability density function (PDF) \( f \in L^2(\mathbb{R}, \mathcal{B}_{\mathbb{R}}, \lambda) \), where \( \lambda \) denotes Lebesgue measure. The cumulative distribution function (CDF) is given by
\begin{equation}
	F(x) = \int_{-\infty}^{x} f(t) \, d\lambda(t), \quad x \in \mathbb{R}.
\end{equation}

We consider a geometric setting where \( X \) denotes the Euclidean distance between two independently and uniformly sampled points on the surface of an \( n \)-dimensional hypersphere of radius \( r > 0 \), denoted by
\begin{equation}
	S_r^n := \left\{ \mathbf{x} \in \mathbb{R}^{n+1} : \|\mathbf{x}\| = r \right\}.
\end{equation}
This random variable represents the chord length between two such points and encodes fundamental aspects of the intrinsic geometry of the hypersphere.

\begin{figure}
	\centering
	\includegraphics[scale=0.5]{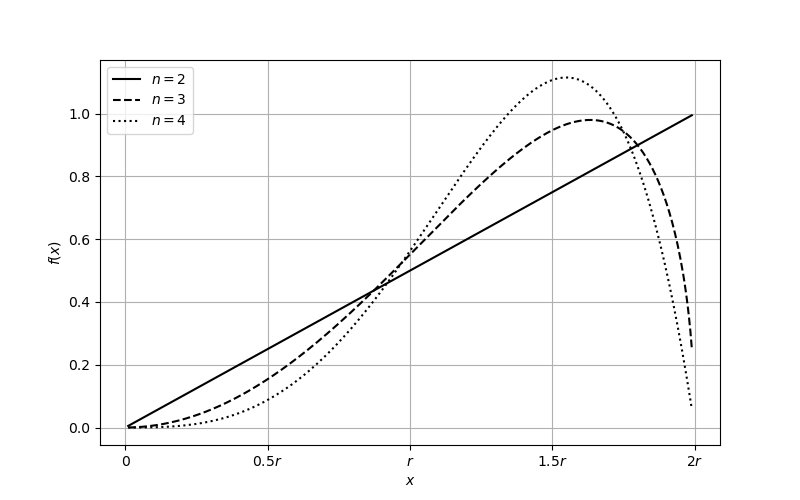}
	\caption{Plots of density functions of the chord length distribution on \( n \)-spheres for various dimensions.}
	\label{Plot_PDF}
\end{figure}

The PDF of the chord length admits the closed-form expression
\begin{equation}
	f(x) = \frac{x}{r^2 \, \mathrm{B}\left( \frac{n}{2}, \frac{1}{2} \right)} \left( \frac{x^2}{r^2} - \frac{x^4}{4r^4} \right)^{\frac{n - 2}{2}}, \quad x \in (0, 2r),
\end{equation}
where \( n \geq 2 \) and \( \mathrm{B}(\cdot, \cdot) \) denotes the Beta function. The support reflects the maximum attainable chord length \( 2r \), corresponding to antipodal points. Figure~\ref{Plot_PDF} depicts the evolution of this distribution as the dimension varies.

The CDF is expressed via the regularized incomplete Beta function as follows:
\begin{equation}
	F(x) = \frac{\mathrm{B}\left( \frac{x^2}{4r^2}; \frac{n}{2}, \frac{n}{2} \right)}{\mathrm{B}\left( \frac{n}{2}, \frac{n}{2} \right)},
\end{equation}
where the incomplete Beta function is defined by
\begin{equation}
	\mathrm{B}(z; a, b) = \frac{1}{\mathrm{B}(a, b)} \int_0^z t^{a-1}(1 - t)^{b - 1} \, dt.
\end{equation}

The moments of \( X \) are computable in closed form. For all \( k > n \), the \( k \)th moment is
\begin{equation}
	\mathbb{E}[X^k] = \frac{2^{k + n - 1}}{\mathrm{B}\left( \frac{n}{2}, \frac{1}{2} \right)} \mathrm{B}\left( \frac{k + n}{2}, \frac{n}{2} \right) r^k.
\end{equation}
In particular, the mean and variance are given by
\begin{equation}
	\mu = \mathbb{E}[X] = \frac{\Gamma^2\left( \frac{n + 1}{2} \right)}{\sqrt{\pi} \, \Gamma\left( n + \frac{1}{2} \right)} \cdot 2^n r,
\end{equation}
\begin{equation}
	\sigma^2 = \operatorname{Var}(X) = \left( 2 - \frac{\Gamma^4\left( \frac{n + 1}{2} \right)}{\pi \, \Gamma^2\left( n + \frac{1}{2} \right)} \cdot 4^n \right) r^2.
\end{equation}

An exceptional feature of the chord length distribution is that its median assumes a universal form, independent of dimension. To identify the median, one solves \( F(m) = \frac{1}{2} \). Substituting \( x = \sqrt{2} \, r \), we find
\begin{equation}
	F(\sqrt{2} \, r) = \frac{\mathrm{B}\left( \frac{1}{2}; \frac{n}{2}, \frac{n}{2} \right)}{\mathrm{B}\left( \frac{n}{2}, \frac{n}{2} \right)} = \frac{1}{2},
\end{equation}
for all \( n \geq 2 \), yielding the closed-form median
\begin{equation}
	m = \sqrt{2} \, r.
\end{equation}

The distribution is strictly unimodal for all \( n \geq 2 \). When \( n = 2 \), the PDF increases monotonically toward \( x = 2r \); for \( n \geq 3 \), the density exhibits a unique interior maximum. Differentiating the PDF yields
\begin{equation}
	f'(x) = \frac{4 r^2 \left[ 4(n - 1) r^2 + (3 - 2n) x^2 \right] \left( \frac{x^2}{r^2} - \frac{x^4}{4 r^4} \right)^{n/2}}{\left(-4 r^2 x + x^3\right)^2 \, \mathrm{B}\left( \frac{n}{2}, \frac{1}{2} \right)}.
\end{equation}
Solving \( f^{\prime}(x) = 0 \) yields the unique mode at
\begin{equation}
	x^* = 2r \sqrt{\frac{n - 1}{2n - 3}}.
\end{equation}
As \( n \to \infty \), this mode converges to \( \sqrt{2} \, r \), aligning with the median. This dimensional convergence illustrates the onset of high-dimensional concentration, a hallmark of hyperspherical geometry.

\subsection{Critical Dimension}

As a result of being a unimodal distribution, the chord length distribution on the \( n \)-sphere permits the application of a classical inequality relating the mean \( \mu \), variance \( \sigma^2 \), and median \( m \) of a unimodal random variable \cite{basu1997mean}:
\begin{equation}
	\frac{(m - \mu)^2}{\sigma^2} \leq \frac{3}{5}.
\end{equation}

Applying this inequality to the chord length distribution yields a nontrivial constraint on the dimensional parameter \( n \). Specifically, it restricts the admissible dimensions for which the statistical spread, as measured by the variance, does not deviate excessively from the central location defined by the median. This provides a dimensionally explicit criterion grounded in the interplay between the shape of the distribution and the geometry of the underlying hypersphere.

\begin{figure}
	\centering
	\includegraphics[scale=0.6]{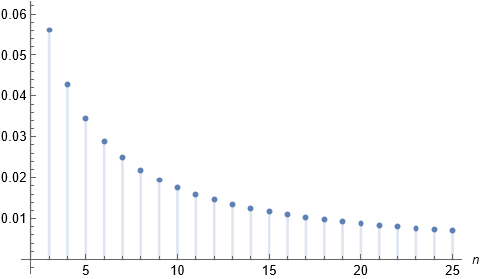}
	\caption{Plot of $\sqrt{2}-C_n.$}
	\label{Plot_Diff}
\end{figure}

Introduce the notation
\begin{equation}
	C_n := \frac{2^n \Gamma^2\left( \frac{n+1}{2} \right)}{\sqrt{\pi} \, \Gamma\left( n + \frac{1}{2} \right)},
\end{equation}
so that \( \mu = C_n \, r \) and \( \sigma^2 = \left( 2 - C_n^2 \right) r^2 \). The squared standardized deviation of the median becomes
\begin{equation}
	\frac{(m - \mu)^2}{\sigma^2} = \frac{(\sqrt{2} - C_n)^2}{2 - C_n^2}.
\end{equation}
By factoring the denominator, the inequality is equivalently rewritten as
\begin{equation}
	\frac{\sqrt{2} - C_n}{\sqrt{2} + C_n} \leq \frac{3}{5},
\end{equation}
which yields the admissible range
\begin{equation}
	\frac{\sqrt{2}}{4} \leq C_n \leq \sqrt{2}.
\end{equation}
We observe that this condition is satisfied for all $n$, but the difference $\sqrt{2}-C_n$, shrinks rapidly with increasing $n$. As the plot of $\sqrt{2}-C_n$ in Figure \ref{Plot_Diff} shows, this difference becomes nearly flat for $n > 19$, signifying a saturation point.

The dimension \( n = 19 \) thus marks a critical threshold in the behavior of the chord length distribution on hyperspheres, signaling the onset of a distinct geometric regime. In lower dimensions, the distribution of interpoint distances retains a rich structural profile, characterized by asymmetry, broad dispersion, and curvature-induced variability. However, as the ambient dimension approaches 19, these features give way to the concentration phenomenon that is emblematic of high-dimensional spaces. At this scale, the geometry of the sphere becomes sufficiently uniform that randomly selected chords cluster tightly around a single dominant length, effectively collapsing the distribution.

This transition reflects a deeper shift from a curvature-governed regime to one dominated by geometric regularity. For \( n \leq 19 \), the statistical properties of interpoint distances are shaped by the global geometry of the sphere, leading to meaningful variation among measures of central tendency. In contrast, for \( n > 19 \), high-dimensional effects suppress this variability: almost all points become nearly equidistant, and the distribution converges sharply toward its mean and median. The identification of this dimensional boundary thus marks the point at which probabilistic diversity yields to geometric determinism, and the sphere's intrinsic uniformity overwhelms random fluctuations.


\subsection{Fisher Information}

We now quantify the amount of information about the radius \( r \) contained in a single observation of the chord length random variable \( X \). This is captured by the Fisher information, which is defined as follows:
\begin{equation}
	\mathcal{I}(r) = \mathbb{E}\left[ \left( \frac{\partial}{\partial r} \log f(X) \right)^2 \right].
\end{equation}
For the chord length distribution on the \( n \)-dimensional hypersphere, this integral admits the closed-form expression
\begin{equation}
	\mathcal{I}(r) = \frac{4(n - 1)n}{(n - 4) r^2}, \quad \text{for } n > 4.
\end{equation}

This expression reveals two structural dependencies: Fisher information decays as \( \mathcal{I}(r) \propto r^{-2} \), and increases roughly linearly with the dimension \( n \) for large $n$. The inverse-square dependence on \( r \) reflects a geometric intuition: as the sphere grows, local perturbations in the radius lead to smaller observable changes in chord lengths, thus weakening the ability to infer \( r \) precisely. Statistically, this implies that the variance of any unbiased estimator of \( r \), constrained by the Cr\'amer-Rao bound, must scale proportionally to \( r^2 \), increasing uncertainty at larger radii.

Conversely, the dimensional growth of \( \mathcal{I}(r) \) embodies the effect of measure concentration. In higher dimensions, interpoint distances become increasingly concentrated around their mean, and the resulting sharp peak in the likelihood function enhances its sensitivity to changes in \( r \). This asymptotic regime is captured by
\begin{equation}
	\mathcal{I}(r) \sim \frac{4n}{r^2}, \quad \text{as } n \to \infty,
\end{equation}
revealing that statistical precision improves unboundedly with dimension, provided the radius is fixed.

\begin{figure}
	\centering
	\includegraphics[scale=0.6]{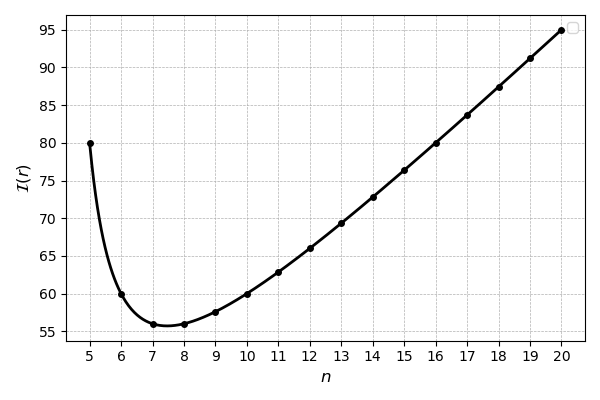}
	\caption{Fisher information as a function of dimension for unit $n$-sphere.}
	\label{Plot_N}
\end{figure}

Interestingly, the Fisher information does not increase monotonically with dimension. As illustrated in Figure \ref{Plot_N}, it attains a global minimum near \( n = 7.46 \). This minimum occurs in close proximity to the integer dimension \( n = 7 \), which we identify as the point of minimal statistical sensitivity for practical purposes. At this dimension, the curvature of the sphere has weakened just enough to suppress observable variability in chord lengths, yet the measure concentration effect has not sufficiently developed to counteract this loss. The resulting dip in Fisher information arises from the confluence of insufficient curvature to generate variability and insufficient dimensionality to enforce concentration, both of which are crucial for enhancing inferential precision.

It is noteworthy that the observed minimum in Fisher information near \( n = 7 \) coincides with a well-known geometric extremum: the volume of the unit \( n \)-sphere achieves its maximum precisely at dimension \( n = 7 \). The volume of an \( n \)-sphere in \( \mathbb{R}^{n+1} \) is given by 
\begin{equation}
	V_n = \dfrac{\pi^{n/2}}{\Gamma(n/2 + 1)},
\end{equation}
which increases with dimension up to \( n = 7 \), after which it decays due to the dominant influence of the Gamma function's growth. This volumetric maximum marks a transition in hyperspherical geometry, from a regime of expansion to one of geometric contraction, and it aligns with the point where the Fisher information is minimized. At this dimension, the hypersphere is most spatially diffuse relative to its embedding space, and hence least sensitive to infinitesimal perturbations in radius. The chord length distribution becomes maximally spread, yet insufficiently concentrated to induce statistical leverage. This geometric-statistical interplay suggests that the Fisher information minimum and volume maximum are not merely coincidental, but reflect a deep duality between curvature, volume, and inferential power on hyperspherical surfaces.

\subsection{Cr\'amer-Rao Inequity}

We now derive the Cr\'amer--Rao lower bound (CRLB) for estimating the radius \( r \) of an \( n \)-dimensional hypersphere. Let \( X_1, \dots, X_m \) be i.i.d.\ observations from the chord length distribution. The Cr\'amer--Rao inequality asserts that the variance of any unbiased estimator \( \hat{r} \) satisfies
\begin{equation}
	\operatorname{Var}(\hat{r}) \geq \frac{1}{m \mathcal{I}(r)},
\end{equation}
where \( \mathcal{I}(r) \) denotes the Fisher information contained in a single observation. For the chord length distribution, this implies that
\begin{equation}
	\operatorname{Var}(\hat{r}) \geq \frac{(n - 4)}{4 m (n - 1)n} r^2, \quad \text{for } n > 4.
\end{equation}

To illustrate the utility of this bound, consider the following sample mean-based estimator
\begin{equation}
	\hat{r} = \frac{\bar{X}}{C_n},
\end{equation}
where
\begin{equation}
	C_n = \frac{2^n \Gamma^2\left( \frac{n+1}{2} \right)}{\sqrt{\pi} \Gamma\left( n + \frac{1}{2} \right)}.
\end{equation}
This estimator is unbiased, since \( \mathbb{E}[X] = C_n r \), which implies \( \mathbb{E}[\hat{r}] = r \). Its variance is given by
\begin{equation}
	\operatorname{Var}(\hat{r}) = \frac{\operatorname{Var}(\bar{X})}{C_n^2} = \frac{\sigma^2}{m C_n^2},
\end{equation}
where \( \sigma^2 = \left( 2 - C_n^2 \right) r^2 \) denotes the variance of the chord length distribution. Substituting this gives
\begin{equation}
	\operatorname{Var}(\hat{r}) = \frac{(2 - C_n^2)}{m C_n^2} r^2.
\end{equation}

Comparing these two expressions confirms that the estimator \( \hat{r} \) is not efficient in general; that is, it does not attain the Cr\'amer--Rao bound. Nonetheless, both the variance and the bound scale linearly with \( r^2 \) and decrease proportionally with the sample size \( m \), in accordance with statistical intuition.

More importantly, as the dimension \( n \to \infty \), we obtain the asymptotic relations
\[
\operatorname{Var}(\hat{r}) \sim \frac{r^2}{2mn}, \quad \text{and} \quad \frac{1}{m \mathcal{I}(r)} \sim \frac{r^2}{4mn}.
\]
This demonstrates that, although suboptimal, the sample mean estimator becomes increasingly efficient in high dimensions. As the dimension \( n \) increases, the chord length distribution concentrates sharply around its limiting value \( \sqrt{2}r \), amplifying Fisher information and reducing variance. Consequently, the Cr\'amer--Rao bound becomes achievable with fewer observations; in fact, the theoretical lower bound on sample size eventually exceeds one, indicating that even a single observation may be sufficient for inference. This striking gain in efficiency is a direct consequence of the geometric regularity induced by high-dimensional hyperspheres.

The resulting sample complexity bound has broader implications for statistical learning in structured spaces. In contrast to the classical expectation that higher dimensions exacerbate overfitting, hyperspherical geometry imposes such strong concentration that observations become nearly deterministic. This introduces an intrinsic form of regularization, independent of data augmentation or penalization, making inference more stable as dimension increases.

These geometric effects naturally benefit machine learning applications that leverage hyperspherical embeddings. In fields such as natural language processing, computer vision, and metric learning, it is increasingly common to constrain latent representations to lie on the unit sphere. The analysis presented here offers theoretical support for the observed robustness of such models: as dimension grows, estimator variance diminishes, enabling high precision from minimal data. Thus, the generalization performance of hyperspherical architectures is not merely empirical, but rooted in geometric compression of variability.

More broadly, these findings challenge traditional views of the curse of dimensionality. While classical capacity measures deteriorate with dimension, the compactness and curvature of high-dimensional spheres can reduce effective statistical complexity, facilitating generalization in expansive spaces. This suggests that ambient curvature may serve as a more appropriate complexity control than traditional dimension-based bounds in manifold-constrained learning.

Ultimately, the hypersphere serves as a prototype for a broader class of curved domains in which geometry itself enhances data efficiency. As manifold-valued models gain prominence, a refined theory of geometry-aware learning will be critical for building scalable and robust inference systems.

\subsection{Characteristic Function}

The characteristic function of the chord length distribution on the \( n \)-sphere of radius \( r \), defined by
\begin{equation}
	\varphi(t) = \mathbb{E}[e^{itX}],
\end{equation}
encapsulates the full probabilistic structure of the random variable $X$, which denotes the Euclidean distance between two independently and uniformly sampled points on the sphere. As the Fourier transform of the probability density function, \( \varphi_n(t) \) encodes moment information, regularity, and tail behavior.

The analytic form of \( \varphi(t) \) depends significantly on whether the ambient dimension \( n \) is even or odd, a consequence of the angular integration structure and the underlying harmonic complexity of hyperspherical geometry.

For \( n = 2 \), the characteristic function admits the closed-form expression
\begin{equation}
	\varphi(t) = \frac{-1 + e^{2irt}(1 - 2irt)}{2r^2 t^2}.
\end{equation}
This rational-exponential form arises from integrating over the circle and reveals limited endpoint smoothness in the corresponding density.

\begin{figure}
	\centering
	\includegraphics[scale=0.45]{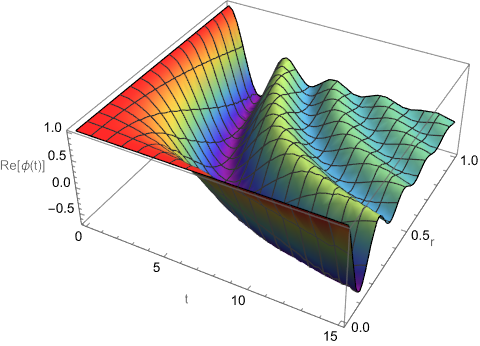}
	\includegraphics[scale=0.45]{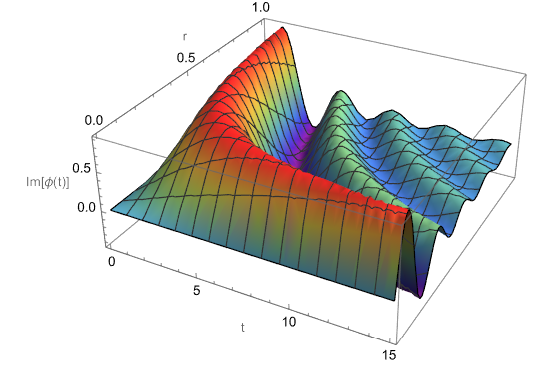}
	\caption{Plots of real (left) and imaginary (right) parts of characteristic function for $n=3$.}
	\label{Plot_Charac}
\end{figure}

In contrast, for \( n = 3 \), the characteristic function is given by
\begin{equation}
	\varphi(t) = \frac{32i r t}{15\pi} + \frac{2}{r^2 t^2} \left[ J_2(2rt) - 2rt\, J_3(2rt) \right] + \frac{2i}{r^2 t^2} \left[ H_2(2rt) - 2rt\, H_3(2rt) \right],
\end{equation}
where \( J_k \) and \( H_k \) denote the Bessel and Struve functions of the first kind. These arise from the Fourier transform of the radially symmetric kernel induced by uniform sampling over the surface, and reflect the spherical harmonic content of odd-dimensional hyperspheres. Their presence signals enhanced smoothness in the density function. Figure \ref{Plot_Charac} depicts the plots of real and imaginary parts of the characteristic function for $n=3$.

The appearance of Bessel and Struve functions in the characteristic function of the chord length distribution on the \( n \)-sphere, particularly in odd dimensions, stems from fundamental harmonic and analytic properties of the underlying geometry. When evaluating the characteristic function, we are effectively computing a Fourier transform of a radially symmetric function over a spherical domain.

In odd dimensions, the Fourier analysis of such radially invariant kernels on the sphere naturally leads to the emergence of special functions. The Bessel functions encode the principal radial oscillations that arise due to the projection of spherical harmonics onto radial bases, while the Struve functions account for boundary behavior and contribute corrections necessary to enforce regularity and normalization. Their joint appearance reflects the precise spectral content of the distribution: the geometry of the sphere imposes angular harmonics whose projections, via Fourier analysis, yield these special functions.

This structure is largely absent in even dimensions, where the angular decomposition leads instead to rational-exponential expressions devoid of harmonic special functions. This distinction reinforces the analytic asymmetry between odd and even-dimensional hyperspheres and reveals how the geometry of the ambient space directly informs the analytic character of the associated statistical distributions.

\newpage

\section{Derivation of Results}
\label{Sec:Results}
\subsection{Probability Density Function}

We first prove that $f(x)$ is a valid probability density function. Positivity of $f(x)$ is immediate from its definition, as all terms are manifestly non-negative over the support $x \in (0, 2r)$. To verify the normalization condition, we compute the integral
\begin{equation}
	\int_0^{2r} f(x)\, dx = \frac{1}{r^2 \, \mathrm{B}\left( \frac{n}{2}, \frac{1}{2} \right)} \int_0^{2r} x \left( \frac{x^2}{r^2} - \frac{x^4}{4r^4} \right)^{\frac{n-2}{2}} \, dx.
\end{equation}
We make the substitution $u = \frac{x^2}{4r^2}$, which implies $x = 2r\sqrt{u}$ and $dx = \frac{r}{\sqrt{u}}\, du$. Substituting, we obtain
\begin{equation}
	\int_0^{2r} f(x)\, dx = \frac{1}{r^2 \, \mathrm{B}\left( \frac{n}{2}, \frac{1}{2} \right)} \int_0^1 (2r \sqrt{u}) \left( 4u - 4u^2 \right)^{\frac{n-2}{2}} \cdot \frac{r}{\sqrt{u}}\, du.
\end{equation}
Simplifying the integrand, we get
\begin{equation}
	\int_0^{2r} f(x)\, dx = \frac{2^{n-1}}{\mathrm{B}\left( \frac{n}{2}, \frac{1}{2} \right)} \int_0^1 u^{\frac{n}{2} - 1}(1 - u)^{\frac{n-2}{2}}\, du.
\end{equation}
This integral is a Beta integral as follows:
\begin{equation}
	\int_0^1 u^{\frac{n}{2} - 1}(1 - u)^{\frac{n-2}{2}}\, du = \mathrm{B}\left( \frac{n}{2}, \frac{n}{2} \right).
\end{equation}
Hence,
\begin{equation}
	\int_0^{2r} f(x)\, dx = \frac{2^{n-1}}{\mathrm{B}\left( \frac{n}{2}, \frac{1}{2} \right)} \cdot \mathrm{B}\left( \frac{n}{2}, \frac{n}{2} \right).
\end{equation}
We now expand both Beta functions using the identity
\begin{equation}
	\mathrm{B}(x, y) = \frac{\Gamma(x)\Gamma(y)}{\Gamma(x+y)}.
\end{equation}
Thus,
\begin{equation}
	\mathrm{B}\left( \frac{n}{2}, \frac{1}{2} \right) = \frac{\Gamma\left( \frac{n}{2} \right)\Gamma\left( \frac{1}{2} \right)}{\Gamma\left( \frac{n+1}{2} \right)} = \frac{\Gamma\left( \frac{n}{2} \right)\sqrt{\pi}}{\Gamma\left( \frac{n+1}{2} \right)},
\end{equation}
and
\begin{equation}
	\mathrm{B}\left( \frac{n}{2}, \frac{n}{2} \right) = \frac{\Gamma^2\left( \frac{n}{2} \right)}{\Gamma(n)}.
\end{equation}
Substituting these into the expression yields
\begin{equation}
	\int_0^{2r} f(x)\, dx = \frac{2^{n-1} \Gamma\left( \frac{n+1}{2} \right)}{\Gamma\left( \frac{n}{2} \right)\sqrt{\pi}} \cdot \frac{\Gamma^2\left( \frac{n}{2} \right)}{\Gamma(n)} = \frac{2^{n-1} \Gamma\left( \frac{n+1}{2} \right) \Gamma\left( \frac{n}{2} \right)}{\sqrt{\pi} \Gamma(n)}.
\end{equation}
We now invoke the duplication formula for the Gamma function given by
\begin{equation}
	\Gamma(n) = \frac{2^{n-1}}{\sqrt{\pi}} \Gamma\left( \frac{n}{2} \right) \Gamma\left( \frac{n+1}{2} \right).
\end{equation}
This shows that the numerator and denominator are equal, hence, implying
\begin{equation}
	\int_0^{2r} f(x)\, dx = 1.
\end{equation}
Thus, $f(x)$ is a valid probability density function.

Next, we derive the chord length distribution from a geometric perspective. Let \( \mathbf{p}_1, \mathbf{p}_2 \in S_r^n \subset \mathbb{R}^{n+1} \) be two points drawn independently and uniformly from the surface of an \( n \)-dimensional hypersphere of radius \( r > 0 \). Our goal is to determine the probability density function of the Euclidean distance \( X = \|\mathbf{p}_1 - \mathbf{p}_2\| \in (0, 2r) \), referred to as the chord length.

Exploiting the rotational symmetry of the hypersphere, we may fix \( \mathbf{p}_1 \) at the north pole,
\begin{equation}
	\mathbf{p}_1 = (0, 0, \ldots, 0, r) \in \mathbb{R}^{n+1},
\end{equation}
and sample \( \mathbf{p}_2 = (p_1, \ldots, p_n, p_{n+1}) \) uniformly from \( S_r^n \). Let \( \theta \in [0, \pi] \) denote the angle between the position vectors \( \mathbf{p}_1 \) and \( \mathbf{p}_2 \), defined via the inner product as follows:
\begin{equation}
	\cos(\theta) = \frac{\mathbf{p}_1 \cdot \mathbf{p}_2}{\|\mathbf{p}_1\| \, \|\mathbf{p}_2\|} = \frac{\mathbf{p}_1 \cdot \mathbf{p}_2}{r^2} = \frac{p_{n+1}}{r}.
\end{equation}
The Euclidean distance is then computed using the identity
\[
\|\mathbf{p}_1 - \mathbf{p}_2\|^2 = \|\mathbf{p}_1\|^2 + \|\mathbf{p}_2\|^2 - 2 \mathbf{p}_1 \cdot \mathbf{p}_2,
\]
which yields
\[
\|\mathbf{p}_1 - \mathbf{p}_2\|^2 = 2r^2(1 - \cos(\theta)),
\]
and hence
\begin{equation}
	X = \|\mathbf{p}_1 - \mathbf{p}_2\| = 2r \sin\left( \frac{\theta}{2} \right).
\end{equation}

To determine the distribution of \( X \), we first derive the density of the angle \( \theta \). The uniform surface measure on \( S_r^n \) implies that the density of \( \theta \in [0, \pi] \), interpreted as the central angle between two random points, is proportional to the surface area of the spherical slice at height \( p_{n+1} = r \cos(\theta) \). This slice corresponds to an \((n-1)\)-dimensional sphere of radius \( \rho = \sqrt{r^2 - p_{n+1}^2} = r \sin(\theta) \).

Since the surface area of an \((n-1)\)-sphere of radius \( \rho \) is given by
\[
A_{n-1}(\rho) = C_{n-1} \rho^{n-1}, \quad C_{n-1} = \frac{2\pi^{n/2}}{\Gamma(n/2)},
\]
the differential surface measure at angle \( \theta \) is proportional to \( \sin^{n-1}(\theta) \, d\theta \). Normalizing over \( [0, \pi] \), we obtain the marginal density of \( \theta \) as follows:
\begin{equation}
	f_\Theta(\theta) = C_n \sin^{n-1}(\theta), \quad \theta \in [0, \pi],
\end{equation}
with normalization constant
\begin{equation}
	C_n = \left( \int_0^\pi \sin^{n-1}(\theta) \, d\theta \right)^{-1} = \frac{\Gamma\left( \frac{n+1}{2} \right)}{\sqrt{\pi} \Gamma\left( \frac{n}{2} \right)}.
\end{equation}

This density reflects the geometric fact that near \( \theta = 0 \) and \( \theta = \pi \), the sphere narrows into caps of small area, whereas near \( \theta = \pi/2 \), the horizontal cross-section reaches its maximal size, resulting in the highest density of angular configurations.

Finally, to compute the probability density function of the chord length \( X \), we apply a change of variables using the transformation
\[
X = 2r \sin\left( \frac{\Theta}{2} \right),
\]
which gives
\begin{equation}
	\theta = 2 \arcsin\left( \frac{x}{2r} \right), \qquad \frac{d\theta}{dx} = \frac{1}{r \sqrt{1 - \left( \frac{x}{2r} \right)^2}}.
\end{equation}
To evaluate the angular density \( f_\Theta(\theta) \) in terms of \( x \), we compute
\begin{equation}
	\sin(\theta) = 2 \sin\left( \frac{\theta}{2} \right) \cos\left( \frac{\theta}{2} \right) = \frac{x}{r} \sqrt{1 - \frac{x^2}{4r^2}},
\end{equation}
so that
\begin{equation}
	\sin^{n-1}(\theta) = \left( \frac{x}{r} \sqrt{1 - \frac{x^2}{4r^2}} \right)^{n-1}.
\end{equation}

Using the standard Jacobian transformation rule, we obtain the density of \( X \) from that of \( \Theta \) as
\begin{equation}
	f_X(x) = f_\Theta(\theta(x)) \left| \frac{d\theta}{dx} \right| = \frac{C_n}{r^n} x^{n-1} \left( 1 - \frac{x^2}{4r^2} \right)^{\frac{n - 2}{2}}, \quad x \in (0, 2r).
\end{equation}

To simplify this expression, we recall the identity
\begin{equation}
	\mathrm{B}\left( \frac{n}{2}, \frac{1}{2} \right) = \frac{\sqrt{\pi} \, \Gamma\left( \frac{n}{2} \right)}{\Gamma\left( \frac{n+1}{2} \right)} \quad \Rightarrow \quad C_n = \frac{1}{\mathrm{B}\left( \frac{n}{2}, \frac{1}{2} \right)},
\end{equation}
and factor the term
\begin{equation}
	x^{n-1} \left( 1 - \frac{x^2}{4r^2} \right)^{\frac{n - 2}{2}} = x \left( \frac{x^2}{r^2} - \frac{x^4}{4r^4} \right)^{\frac{n - 2}{2}} r^{n - 2}.
\end{equation}
Substituting back, we arrive at the canonical form of the chord length distribution as follows:
\begin{equation}
	f(x) = \frac{x}{r^2 \, \mathrm{B}\left( \frac{n}{2}, \frac{1}{2} \right)} \left( \frac{x^2}{r^2} - \frac{x^4}{4r^4} \right)^{\frac{n - 2}{2}}, \qquad x \in (0, 2r).
\end{equation}
This completes the derivation.

\subsection{Cumulative Distribution Function}

To compute the cumulative distribution function $F(x)$, we integrate
\begin{equation}
	F(x) = \int_0^x f(u) \, du.
\end{equation}
We make the change of variables $z = \frac{u^2}{4r^2}$, so that $u = 2r \sqrt{z}$ and $du = \frac{r}{\sqrt{z}} \, dz$. Note that
\begin{equation}
	\frac{u^2}{r^2} - \frac{u^4}{4r^4} = \frac{4r^2 z - 16r^2 z^2}{r^2} = 4z(1 - z),
\end{equation}
and $u = 2r \sqrt{z}$, so
\begin{equation}
	F(x) = \frac{1}{r^2 \, \mathrm{B}\left( \frac{n}{2}, \frac{1}{2} \right)} \int_0^{x} u \left( \frac{u^2}{r^2} - \frac{u^4}{4r^4} \right)^{\frac{n - 2}{2}} \, du.
\end{equation}
Substituting in the new variable gives
\begin{equation}
	F(x) = \frac{1}{r^2 \, \mathrm{B}\left( \frac{n}{2}, \frac{1}{2} \right)} \int_0^{\frac{x^2}{4r^2}} (2r \sqrt{z}) \left( 4z(1 - z) \right)^{\frac{n - 2}{2}} \cdot \frac{r}{\sqrt{z}} \, dz.
\end{equation}
Simplifying, we get
\begin{equation}
	F(x) = \frac{2^{n-1}}{\mathrm{B}\left( \frac{n}{2}, \frac{1}{2} \right)} \int_0^{\frac{x^2}{4r^2}} z^{\frac{n}{2} - 1} (1 - z)^{\frac{n - 2}{2}} \, dz.
\end{equation}
This integral is the incomplete Beta function, yielding
\begin{equation}
	F(x) = \frac{2^{n-1}}{B\left(\frac{n}{2}, \frac{1}{2}\right)} \mathrm{B}\left( \frac{x^2}{4r^2}; \frac{n}{2}, \frac{n}{2} \right),
\end{equation}
which can be equivalently expressed in terms of the regularized incomplete Beta function as follows:
\begin{equation}
	F(x) = \frac{\mathrm{B}\left( \frac{x^2}{4r^2}; \frac{n}{2}, \frac{n}{2} \right)}{\mathrm{B}\left( \frac{n}{2}, \frac{n}{2} \right)}.
\end{equation}

\subsection{Raw Moments}
To compute the $k$th raw moment of $X$, we evaluate
\begin{equation}
	\mathbb{E}[X^k] = \int_0^{2r} x^k f(x) \, dx = \frac{1}{r^2 \, \mathrm{B}\left( \frac{n}{2}, \frac{1}{2} \right)} \int_0^{2r} x^{k+1} \left( \frac{x^2}{r^2} - \frac{x^4}{4r^4} \right)^{\frac{n - 2}{2}} \, dx.
\end{equation}
Let $u = \frac{x^2}{4r^2}$, so that $x = 2r \sqrt{u}$ and $dx = \frac{r}{\sqrt{u}} \, du$. The substitution yields
\begin{equation}
	\mathbb{E}[X^k] = \frac{1}{r^2 \, \mathrm{B}\left( \frac{n}{2}, \frac{1}{2} \right)} \int_0^1 (2r \sqrt{u})^{k+1} (4u - 4u^2)^{\frac{n - 2}{2}} \cdot \frac{r}{\sqrt{u}} \, du.
\end{equation}
Simplifying the integrand we get
\begin{equation}
	\mathbb{E}[X^k] = \frac{2^{k+n-1} r^k}{\mathrm{B}\left( \frac{n}{2}, \frac{1}{2} \right)} \int_0^1 u^{\frac{k+n}{2} - 1} (1 - u)^{\frac{n - 2}{2}} \, du.
\end{equation}
The integral is recognized as the following Beta function
\begin{equation}
	\int_0^1 u^{\frac{k+n}{2} - 1} (1 - u)^{\frac{n - 2}{2}} \, du = \mathrm{B}\left( \frac{k+n}{2}, \frac{n}{2} \right).
\end{equation}
Therefore, the moment is given by
\begin{equation}
	\mathbb{E}[X^k] = \frac{2^{k+n-1}}{\mathrm{B}\left( \frac{n}{2}, \frac{1}{2} \right)} \mathrm{B}\left( \frac{k+n}{2}, \frac{n}{2} \right) r^k.
\end{equation}

\section{Conclusion} \label{Sec:Conclusion}

We have analyzed the statistical geometry of chord lengths on \( n \)-dimensional hyperspheres through explicit derivations of their distributional, inferential, and harmonic properties. Beginning with a closed-form expression for the probability density function, we explored its cumulative distribution, moments, and unimodal structure, emphasizing how these features evolve with dimension.

A critical transition was identified at dimension \( n = 19 \), marking the onset of high-dimensional concentration where interpoint distances become nearly deterministic. This regime shift was further quantified through the behavior of the Fisher information, which was shown to vary non-monotonically with dimension and to attain a minimum near \( n = 7 \), coinciding with the maximum volume of the unit sphere. This alignment underscores a duality between geometric expansion and inferential sensitivity on curved domains.

Our analysis also revealed an analytic dichotomy in the characteristic function of the chord length distribution: in even dimensions, the function assumes a rational-exponential form, while in odd dimensions it involves Bessel and Struve functions. This reflects deeper harmonic asymmetries between different dimensional regimes.

The results presented here offer a unified framework for understanding how curvature and dimension regulate statistical behavior on hyperspherical surfaces. Beyond their theoretical interest, these findings provide insight into the geometric structure of data in high-dimensional settings, with potential implications for inference, learning, and representation in manifold-constrained models.

Future work may extend these ideas to more general Riemannian manifolds, incorporate constraints such as geodesic restrictions or caps, or explore inverse problems such as estimating dimension or curvature from observed interpoint distances.


\begin{thebibliography}{1}
	
	\bibitem{basu1997mean}
	Basu, Sanjib and DasGupta, Anirban.
	\newblock The mean, median, and mode of unimodal distributions: a characterization.
	\newblock {\em Theory of Probability and Its Applications}, 41(2):210--223, 1997. SIAM.
	
\end{thebibliography}

\end{document}